# Non-Commutative Simplicial Complexes and the Baum-Connes Conjecture

Joachim Cuntz
Mathematisches Institut
Universität Münster
Einsteinstr. 62
48149 Münster
cuntz@math.uni-muenster.de

September 8, 2018

**Abstract**

We associate a non-commutative $C^*$-algebra with any locally finite simplicial complex. We determine the $K$-theory of these algebras and show that they can be used to obtain a conceptual explanation for the Baum-Connes conjecture.

## 1 Introduction

The Baum-Connes conjecture predicts a formula for the K-theory of group $C^*$-algebras $C^*_{red}\Gamma$ and, more generally, for the $K$-theory of reduced crossed product $C^*$-algebras $A \rtimes_r \Gamma$, [2], [3], [4], [5]. According to the formula, this $K$-theory should be isomorphic, via a map defined by Baum-Connes, to a topologically defined $K$-theory, with coefficients in $A$, of a classifying space $\underline{E}_\Gamma$ for proper actions of the given group.

The Baum-Connes conjecture contains the Novikov conjecture and the generalized Kadison conjecture and plays an important motivating role in current research on topological $K$-theory. While counterexamples to the conjecture have recently been constructed (at least for crossed products), by various authors (Higson, Lafforgue-Skandalis, Yu), it is still expected to be true for a very large class of groups.

The conjecture had been elaborated by Baum and Connes over the years with amazing insight, based mainly on the evidence given by a certain (at the time quite limited) number of examples. It is remarkable that it has indeed been verified for large and important classes of groups (see e.g. [10], [7], [11]). There was however no compelling explanation for the choice of the left hand side of the conjecture, i.e. the topological $K$-theory of $\underline{E}_\Gamma$.

---

*Research supported by the Deutsche Forschungsgemeinschaft

In this article we present such a conceptual explanation. We analyze projections and in particular projections with only finitely many non-zero coefficients in crossed product algebras, or more generally in algebras with a $\Gamma$-grading. We then show that the coefficients of such projections lead naturally to representations of a noncommutative classifying space which is exactly the noncommutative analogue of $\underline{E}_\Gamma$. More precisely, we obtain a projective system of $C^*$-algebras which are a very straightforward noncommutative generalization of the algebra of functions on certain simplicial complexes. Such a $C^*$-algebra can be associated naturally with any locally finite simplicial complex.

We develop to some extent the theory of these noncommutative simplicial complexes which, we think, have an interest of their own. They do seem to be very natural objects in connection with noncommutative topology. We mention that similar noncommutative simplicial complexes have been used by W.Winter, [13] to obtain a definition of covering dimension for nuclear $C^*$-algebras that has good properties. We show in particular that, under certain finiteness conditions, the noncommutative algebras associated with a simplicial complex have the same $K$-theory as their classical counterpart. This can also be done equivariantly - a fact, which is important in connection with the Baum-Connes conjecture.

In section 3 we investigate another variant of a noncommutative $C^*$-algebra that can be associated with a simplicial complex. As we show later, these algebras, too, have a natural connection to a version of finite $K$-theory for algebras with a $\Gamma$-grading. Their $K$-theory however has less good properties. We compute the $K$-theory for a natural example.

I am grateful to Georges Skandalis for hospitality and very valuable discussions during my stay at Paris. I am also indebted to Siegfried Echterhoff and Ralf Meyer for important comments. The research for this article has been supported by the DFG through the SFB 478 and partially also by the European Research Training Network HPRN-CT-1999-00118 .

## 2 Noncommutative simplicial complexes.

Recall that a simplicial complex $\Sigma$ can be defined as a set of non-empty finite subsets (the simplexes in $\Sigma$) of a set $V_\Sigma$ (the set of vertices of $\Sigma$) satisfying the following axioms:

- if $s \in V_\Sigma$ then $\{s\} \in \Sigma$
- if $F \in \Sigma$ and $\emptyset \neq E \subset F$, then $E \in \Sigma$.

$\Sigma$ is called locally finite, if every vertex of $\Sigma$ is contained in only finitely many simplexes of $\Sigma$. In the following, we will always assume that $\Sigma$ is a locally finite simplicial complex. With a simplicial complex $\Sigma$ one can associate canonically a topological space $|\Sigma|$ - the "geometric realization" of the complex. It can be defined as the space of maps $f : V_\Sigma \to [0,1]$ with finite support such that

$$\sum_{s \in V_\Sigma} f(s) = 1$$

with the topology of pointwise convergence. If $\Sigma$ is locally finite, then $|\Sigma|$ is locally compact.

With a locally finite complex we associate in this section two canonical $C^*$-algebras:

(1) $\mathcal{C}_\Sigma$ is the universal $C^*$-algebra with positive generators $h_s, s \in V_\Sigma$ satisfying the relations

  - $h_{a_0} h_{a_1} \ldots h_{a_n} = 0$ whenever $\{a_0, \ldots, a_n\}$ is not a simplex in $\Sigma$ (note that here repetitions in the $a_i$ are allowed).
  - $\sum_{s \in V_\Sigma} h_s h_t = h_t$, $t \in V_\Sigma$

  The sum in the second condition is necessarily finite and this condition may be abbreviated formally to $\sum_{s \in V_\Sigma} h_s = 1$.

(2) $\mathcal{C}_\Sigma^{ab}$ is the universal $C^*$-algebra with generators and relations as in (1) but satisfying in addition $h_s h_t = h_t h_s$ for all $s, t \in V_\Sigma$.

Homomorphisms from $\mathcal{C}_\Sigma$ (or $\mathcal{C}_\Sigma^{ab}$) are easily defined by specifying the images of the generators. Whenever $h'_s, s \in V_\Sigma$ are positive elements in some $C^*$-algebra $A$ satisfying the relations in (1) (or in (2)), then there is a unique *-homomorphism from $\mathcal{C}_\Sigma$ (or from $\mathcal{C}_\Sigma^{ab}$) to $A$ mapping each $h_s$ to $h'_s$.

The algebra $\mathcal{C}_\Sigma^{ab}$ is just the abelianization of $\mathcal{C}_\Sigma$ and we have a canonical surjective map

$$\mathcal{C}_\Sigma \longrightarrow \mathcal{C}_\Sigma^{ab}$$

The natural maps between simplicial complexes $\Sigma$ and $\Sigma'$ in our context are the proper simplicial maps. Such a map is given by a map $\varphi : V_\Sigma \to V_{\Sigma'}$ such that $\{\varphi(s_0), \varphi(s_1), \ldots, \varphi(s_n)\}$ is a simplex in $\Sigma'$, whenever $\{s_0, s_1, \ldots, s_n\}$ is a simplex in $\Sigma$ and such that the preimage of any $s \in V_{\Sigma'}$ is finite. Any such map between locally finite simplicial complexes induces a homomorphism $\mathcal{C}_{\Sigma'} \to \mathcal{C}_\Sigma$ mapping a generator $h_s$ to $\sum_{\varphi(t)=s} h_t$ ($h_s$ is mapped to 0 if $s$ is not in the image of $\varphi$). Important special cases, which will be used extensively, are the evaluation maps which are induced by the inclusion of a subcomplex. $\Sigma_0$ is a subcomplex of $\Sigma$ if $V_{\Sigma_0} \subset V_\Sigma$ and any simplex in $\Sigma_0$ is also a simplex in $\Sigma$. There is a natural evaluation homomorphism

$$\pi_{\Sigma_0, \Sigma} : \mathcal{C}_\Sigma \longrightarrow \mathcal{C}_{\Sigma_0}$$

which sends a generator $h_s, s \in V_{\Sigma_0}$ of $\mathcal{C}_\Sigma$ to the corresponding generator $h_s$ of $\mathcal{C}_{\Sigma_0}$ and sends the other generators $h_s, s \notin V_{\Sigma_0}$ to 0.

**Proposition 2.1** *The (abelian) $C^*$-algebra $\mathcal{C}_\Sigma^{ab}$ is isomorphic to the algebra $\mathcal{C}_0(|\Sigma|)$ of continuous functions vanishing at infinity on the geometric realization $|\Sigma|$ of $\Sigma$.*

*Proof.* We know that $\mathcal{C}_\Sigma^{ab} \cong \mathcal{C}_0(X)$ where $X = \operatorname{Spec} \mathcal{C}_\Sigma^{ab}$. By definition, the elements of $X$ are the homomorphisms $\mathcal{C}_\Sigma^{ab} \to \mathbb{C}$. Such a homomorphism is uniquely determined by the values it takes on the generators $h_s$. These values are real numbers $t_s$ between 0 and 1 with sum equal to 1. The set of all $s$ for which $t_s$ is non-zero form necessarily a simplex in $\Sigma$, since then the product of the corresponding $h_s$ has to be non-zero. Thus we get a bijection between such homomorphisms and maps $f: V_\Sigma \to [0,1]$ supported on a simplex in $\Sigma$ satisfying $\sum_{s \in V_\Sigma} f(s) = 1$. This is exactly the description of the geometric realization for $\Sigma$ given above. $\square$

The $C^*$-algebras $\mathcal{C}_\Sigma$ and $\mathcal{C}_\Sigma^{ab}$ have a natural filtration by ideals which reduces to the skeleton filtration in the abelian case. For each $n$, we let $I_n$ denote the ideal in $\mathcal{C}_\Sigma$ generated by all products of generators $h_s$ containing at least $n+1$ pairwise different generators. We also denote by $I_n^{ab}$ the image of $I_n$ in $\mathcal{C}_\Sigma^{ab}$.

The subcomplex $\Sigma^{n-1}$ of $\Sigma$ that has the same vertices and whose simplexes are exactly the $k$-simplexes of $\Sigma$ with $k \leq n-1$ is called the $n-1$-skeleton of $\Sigma$.

**Proposition 2.2** *The ideal $I_n$ is exactly the kernel of the evaluation map*

$$\pi_{\Sigma^{n-1},\Sigma}: \mathcal{C}_\Sigma \longrightarrow \mathcal{C}_{\Sigma^{n-1}}$$

*The abelian $C^*$-algebra $\mathcal{C}_\Sigma^{ab}/I_n^{ab}$ is isomorphic to $\mathcal{C}_0(|\Sigma^{n-1}|)$.*

*Proof.* Obvious from the proof of 2.1. $\square$

Let $\Delta$ be an $n$-simplex with vertices $\{t_0, \ldots, t_n\}$ and $\Delta^i, i = 0, \ldots, n$ the $(n-1)$-face of $\Delta$ with vertices $\{t_0, \ldots, t_n\} \setminus \{t_i\}$. By definition, $\mathcal{C}_\Delta$ is the unital $C^*$-algebra with generators $h_s$, $s \in \{t_0, \ldots, t_n\}$ such that $h_s \geq 0$ and $\sum h_s = 1$. There are natural evaluation maps $\pi_{\Delta^i,\Delta}: \mathcal{C}_\Delta \to \mathcal{C}_{\Delta^i}$ mapping the generator $h_{t_i}$ to 0. We denote by $J_\Delta$ the ideal in $\mathcal{C}_\Delta$ generated by products of generators containing all the $h_{t_i}, i = 0, \ldots, n$. Note that $J_\Delta$ is exactly the intersection of the kernels of all the $\pi_{\Delta^i,\Delta}$.

This is nothing but the ideal $I_n$ in $\mathcal{C}_\Delta$. The abelian version $J_\Delta^{ab}$ is obviously isomorphic to the algebra of continuous functions on the Euclidean $n$-simplex

$$|\Delta| \cong \{(x_0, \ldots, x_n) \in \mathbb{R}^{n+1} \mid 0 \leq x_i \leq 1, \sum_{i=0}^n x_i = 1\}$$

that vanish on the boundary.

**Lemma 2.3** *The ideal $J_\Delta$ is essential in $\mathcal{C}_\Delta$.*

*Proof.* We proceed by induction on $n$ and assume that the assertion has been proved already for $n-1$ (it is obviously true for $n = 0, 1$). Given $t_i$ in the vertex set $\{t_0, \ldots, t_n\}$ for $\Delta$, and $t \in [0,1]$, we define a homomorphism $\varphi_t: \mathcal{C}_\Delta \to \mathcal{C}_{\Delta^i}$ by $\varphi_t(h_s) = (1-t)h_s, s \neq t_i$, and $\varphi_t(h_{t_i}) = t1$. Let $\varphi$ denote the corresponding homomorphism from $\mathcal{C}_\Delta$ into the algebra $\mathcal{C}_{\Delta^i}[0,1]$ of continuous functions on the unit interval with values in $\mathcal{C}_{\Delta^i}$. The image of $J_\Delta$ under $\varphi$ equals $J_{\Delta^i}(0,1)$, i.e. the functions with values in $J_{\Delta^i}$ vanishing in 0 and 1. If $xJ_\Delta = \{0\}$ for $x \in \mathcal{C}_\Delta$, then $\varphi(x)J_{\Delta^i}(0,1) = \{0\}$. Since $J_{\Delta^i}$ is essential in $\mathcal{C}_{\Delta^i}$ by hypothesis, this implies that $\varphi_t(x) = 0$ for $0 < t < 1$ and by continuity also that $\varphi_0(x) = \pi_{\Delta^i,\Delta}(x) = 0$. It follows that $x \in J_\Delta = \bigcap_i \operatorname{Ker} \pi_{\Delta^i,\Delta}$ whence $x = 0$. □

As a consequence of the previous lemma, the natural map from $\mathcal{C}_\Delta$ into the multiplier algebra $\mathcal{M}(J_\Delta)$ is injective.

**Lemma 2.4** *Let $\Sigma$ be a locally finite simplicial complex and $I_n$ be the ideal in $\mathcal{C}_\Sigma$ defined above. Then*

$$I_n/I_{n+1} \cong \bigoplus_\Delta J_\Delta$$

*where $J_\Delta \subset \mathcal{C}_\Delta$ is the ideal defined above and the direct sum is taken over all $n$-simplexes $\Delta$ in $\Sigma$. The natural map $\mathcal{C}_\Sigma/I_{n+1} \to \mathcal{M}(J_\Delta)$, for each $\Delta$, factors as*

$$\mathcal{C}_\Sigma/I_{n+1} \xrightarrow{\pi_{\Delta,\Sigma^n}} \mathcal{C}_\Delta \underset{\cong}{\subseteq} \mathcal{M}(J_\Delta)$$

*Proof.* The first part of the assertion follows from the fact that, mod $I_{n+1}$, any product containing $n+1$ different generators is orthogonal to any further generator. The second part follows from Lemma 2.3 and from the fact that the natural map into the multiplier algebra of $J_\Delta$ annihilates all $h_s$ for $s \notin \Delta$. □

Call a simplex of $\Sigma$ maximal if it is not strictly contained in any other simplex of $\Sigma$. Obviously, in a locally finite complex, any simplex is contained in a maximal simplex.

**Lemma 2.5** *Let $\Sigma$ be a locally finite simplicial complex and $x \in \mathcal{C}_\Sigma$. If $\pi_{\Delta,\Sigma}(x) = 0$ for all simplexes $\Delta$ in $\Sigma$, then $x = 0$.*
*The same holds if $\pi_{\Delta,\Sigma}(x) = 0$ for all maximal simplexes $\Delta$ in $\Sigma$.*

*Proof.* It follows from 2.4 by induction on $n$ that $x$ has to lie in $I_n$ for each $n$. On the other hand, for a locally finite complex, one has for each generator $h_s$, that $h_s I_k = 0$ for sufficiently large $k$. Thus, if $z$ is any polynomial in the generators of $\mathcal{C}_\Sigma$ such that $\|x - z\| < \varepsilon$, then $\overline{\mathcal{C}_\Sigma z \mathcal{C}_\Sigma} \cap I_k = 0$ for sufficiently large $k$. This implies that the image of $z$ in the quotient by $I_k$ has the same norm as $z$, so that $\|x\| < \varepsilon$ for any $\varepsilon$.

The second statement follows from the fact that, if $F$ is a face of $\Delta$, then $\pi_{F,\Sigma} = \pi_{F,\Delta} \circ \pi_{\Delta,\Sigma}$. □

**Remark 2.6** *As a consequence of the preceding lemma, we obtain, for any $x$ in $\mathcal{C}_\Sigma$*

$$\|x\| = \sup\{ \|\pi_{\Delta,\Sigma}(x)\| \mid \Delta \text{ is a maximal simplex in } \Sigma\}$$

It follows that, for any locally finite simplicial complex $\Sigma$, the algebra $\mathcal{C}_\Sigma$ is a fibered product, over all simplexes $\Delta$ in $\Sigma$, of algebras $\mathcal{C}_\Delta$.

**Theorem 2.7** *The algebra $\mathcal{C}_\Sigma$ is isomorphic to the subalgebra $A_\Sigma$ of the direct sum $\bigoplus \mathcal{C}_\Delta$ over all simplexes in $\Sigma$, consisting of families $(x_\Delta)_{\Delta \in \Sigma}$ such that*

$$\pi_{F,\Delta}(x_\Delta) = \pi_{F,\Delta'}(x_{\Delta'})$$

*whenever $F$ is a common face of $\Delta$ and $\Delta'$.*
*In the direct sum it suffices to take the sum over all maximal simplexes.*

*Proof.* It follows from Lemma 2.5 that the map $\sigma = \prod \pi_{\Delta,\Sigma} : \mathcal{C}_\Sigma \to \bigoplus \mathcal{C}_\Delta$ is injective. Its image is obviously in $A_\Sigma$.
On the other hand, if $x = (x_\Delta)$ is an element of $A_\Sigma$ and $\Delta$ is a maximal simplex in $\Sigma$, we can find $y_\Delta$ in the subalgebra of $\mathcal{C}_\Sigma$ generated by the $h_s$, $s \in \Delta$, such that $x_\Delta = \pi_{\Delta,\Sigma}(y_\Delta)$. Similarly, an element of $\mathcal{C}_\Delta$, that vanishes under the evaluation on a face $F$, can be lifted to an element in the ideal generated by the $h_s$, $s \notin F$ in the subalgebra of $\mathcal{C}_\Sigma$ generated by the $h_s$, $s \in \Delta$. By iterating this, we can find $y$ in $\mathcal{C}_\Sigma$ such that $x_\Delta = \pi_{\Delta,\Sigma}(y)$ for any finite collection of maximal simplexes $\Delta$.
It follows that, for any $\varepsilon > 0$, there is $y \in \mathcal{C}_\Sigma$ such that $\|x - \sigma(y)\| < \varepsilon$. $\square$

**Lemma 2.8** *Let $\Sigma$ be a locally finite simplicial complex and $\Delta$ a simplex in $\Sigma$. The canonical evaluation map $\pi_\Delta : \mathcal{C}_\Sigma \to \mathcal{C}_\Delta$ admits a canonical completely positive lift $\varphi_\Delta : \mathcal{C}_\Delta \to \mathcal{C}_\Sigma$. This lift is equivariant for the action of any group on $\Sigma$ mapping $\Delta$ into itself.*

*Proof.* Let $V$ be the set of vertices of $\Sigma$, $h_s, s \in V$ the canonical generators of $\mathcal{C}_\Sigma$ and $\Delta = \{s_0, s_1, \ldots, s_n\}$.
Let $D$ be the $C^*$-subalgebra of $\mathcal{C}_\Sigma$ generated by $h_{s_0}, \ldots, h_{s_n}$ and denote by $h'_0, h'_1, \ldots, h'_n$ the generators of $\mathcal{C}_\Delta$. We define a homomorphism $\alpha : \mathcal{C}_\Delta \to \tilde{D}$ by putting

$$\alpha(h'_i) = h_{s_i} + \tfrac{1}{n+1}(1 - \sum_{i=0}^{n} h_{s_i})$$

The desired map $\varphi_\Delta$ can now be obtained as

$$\varphi_\Delta(x) = (h_{s_0} + \ldots + h_{s_n})\alpha(x)(h_{s_o} + \ldots + h_{s_n}) \in D \subset \mathcal{C}_\Delta.$$

$\square$

**Lemma 2.9** *Let $\Delta = \{s_0, s_1, \ldots, s_n\}$ be a simplex and $F^i = \{s_0, s_1, \ldots, s_n\} \setminus \{s_i\}$, $i = 0, 1, \ldots, n$ its faces. Let $K_i$ be the kernel of the evaluation map $\mathcal{C}_\Delta \to \mathcal{C}_{F^i}$, i.e. the ideal generated by $h_i$.*
*Then there are multipliers $M_i$, $i = 0, 1, \ldots, n$, of $J_\Delta$ such that $0 \leq M_i \leq 1$, $\sum M_i = 1$, $M_i$ commute with $\mathcal{C}_\Delta$ modulo $J_\Delta$ and $M_i K_j \subset J_\Delta$ for $j \neq i$.*
*Moreover, the $M_i$ can be chosen equivariant for the action of any group on $\Delta$.*

*Proof.* This is exactly Kasparov's technical theorem, [8]. □

**Proposition 2.10** *Let $\Gamma$ be a group acting on the locally finite complex $\Sigma$ by simplicial maps. The quotient map $\mathcal{C}_\Sigma/I_{n+1} \to \mathcal{C}_\Sigma/I_n$ admits an equivariant completely positive splitting.*

*Proof.* The map $\mathcal{C}_\Sigma/I_{n+1} \to \mathcal{C}_\Sigma/I_n$ corresponds to the inclusion map $\Sigma^{n-1} \to \Sigma^n$ between the skeletons.
By Theorem 2.7, $\mathcal{C}_{\Sigma^{n-1}}$ is isomorphic to a subalgebra $A_{\Sigma^{n-1}}$ of the direct sum $\bigoplus_{\Delta' \in E'} \mathcal{C}_{\Delta'}$, where $E'$ is the set of all maximal simplexes $\Delta'$ in $\Sigma^{n-1}$. Similarly, $\mathcal{C}_{\Sigma^n}$ is contained in the direct sum $\bigoplus_{\Delta \in E} \mathcal{C}_\Delta$, where $E$ is the set of all maximal simplexes $\Delta$ in $\Sigma^n$. Note that some of the simplex algebras appear in both direct sums.
We define a completely positive map
$$\psi : \bigoplus \mathcal{C}_{\Delta'} \longrightarrow \bigoplus \mathcal{C}_\Delta$$
as a sum over maps $\psi_{\Delta,\Delta'} : \mathcal{C}_{\Delta'} \to \mathcal{C}_\Delta$, where, given $\Delta'$ in $E'$, $\Delta$ is in the set $M_{\Delta'}$ of all simplexes $\Delta$ in $E$ such that $\Delta \supset \Delta'$, and where

$$\psi_{\Delta',\Delta'} = \mathrm{id} \quad \text{if } \Delta' \text{ is not a face of an } n\text{-simplex in } \Sigma^n \text{ i.e. if } M_{\Delta'} = \{\Delta'\}$$

$$\psi_{\Delta,\Delta'}(x) = M_F^{\frac{1}{2}} \varphi_F(x) M_F^{\frac{1}{2}} \quad \text{if } \Delta' = F \text{ is a face of the } n\text{-simplex } \Delta \text{ in } \Sigma^n$$

Here, $\varphi_F$ is as in 2.8 and $M_F$ is as in 2.9. It is easily checked that the sum $\psi$ over the $\psi_{\Delta,\Delta'}$ makes sense and that $\psi$ is a lift for the quotient map $\mathcal{C}_{\Sigma^n} \to \mathcal{C}_{\Sigma^{n-1}}$. □

An equivariant homomorphism $A_1 \to A_2$ between two $C^*$-algebras equipped with an action of the group $G$ induces a $KK^G$-equivalence, if, for any $G$-algebra $B$ the induced maps $KK^G(A_2, B) \to KK^G(A_1, B)$ and $KK^G(B, A_1) \to KK^G(B, A_2)$ are isomorphisms.

**Theorem 2.11** *Let $G$ be a group acting by simplicial maps on a simplex $\Delta$ and on a finite simplicial complex $\Sigma$.*

(a) *There is an equivariant map $\alpha : \mathcal{C}_\Delta \to \mathbb{C}$ which is an equivariant homotopy inverse to the natural inclusion $\mathbb{C} \to \mathcal{C}_\Delta$ and thus induces a $KK^G$-equivalence.*

(b) *The canonical map $J_\Delta \to J_\Delta^{ab}$ induces a $KK^G$-equivalence.*

(c) *The canonical map $\mathcal{C}_\Sigma \to \mathcal{C}_\Sigma^{ab}$ induces a $KK^G$-equivalence.*

*Proof.* (a) If $h_0, h_1, \ldots, h_n$ denote the generators of $\mathcal{C}_\Delta$, we put $\alpha(h_i) = \frac{1}{n+1} 1$. If we compose $\alpha$ with the inclusion $\mathbb{C} \to \mathcal{C}_\Delta$ it becomes homotopic to the identity map via the equivariant homotopy defined by

$$\alpha_t(h_i) = th_i + \tfrac{1-t}{n+1} 1, \qquad t \in [0,1]$$

We prove (b) and (c) simultaneously by induction on the dimension $n$ of $\Delta$ and $\Sigma$. Both assertions are trivially true for $n = 0$.

Let $\Delta$ be an $n$-dimensional simplex and $S$ its boundary, i.e. the standard simplicial $n-1$-sphere. There is an extension
$$0 \longrightarrow J_\Delta \longrightarrow \mathcal{C}_\Delta \longrightarrow \mathcal{C}_S \longrightarrow 0$$
and the analogous extension for the abelianized algebras. By (a) the map $\mathcal{C}_\Delta \to \mathcal{C}_\Delta^{ab}$ is a $KK^G$-equivalence and by induction hypothesis, the map $\mathcal{C}_S \to \mathcal{C}_S^{ab}$ is so, too. Therefore, by the five lemma and using 2.10, the third map $J_\Delta \to J_\Delta^{ab}$ is a $KK^G$-equivalence.

Let now $(I_k)$ be the skeleton filtration for the $n$-dimensional finite complex $\mathcal{C}_\Sigma$. We have $I_{n+1} = 0$ and $I_n \cong \bigoplus J_\Delta$ where the sum is over all $n$-dimensional simplexes in $\Sigma$. There is an extension
$$0 \longrightarrow I_n \longrightarrow \mathcal{C}_\Sigma \longrightarrow \mathcal{C}_{\Sigma^{n-1}} \longrightarrow 0$$
where $\Sigma^{n-1}$ is the $n-1$-skeleton of $\Sigma$. Since the assertion is true for $\mathcal{C}_{\Sigma^{n-1}}$, by induction hypothesis, and for $I_n$, by (b), which we just proved, we can again use the five lemma, to conclude that the map is a $KK^G$-equivalence. $\square$

**Proposition 2.12** *Let $\Sigma$ be a locally finite simplicial complex and $(I_k)$ the skeleton filtration for $\mathcal{C}_\Sigma$. Assume that $G$ is a group acting on $\Sigma$ (by simplicial maps). Then the map $I_n/I_{n+1} \to I_n^{ab}/I_{n+1}^{ab}$ is a $KK^G$-equivalence for each $n$.*

*Proof.* We know that
$$I_n/I_{n+1} \cong \bigoplus_\Delta J_\Delta \qquad I_n^{ab}/I_{n+1}^{ab} \cong \bigoplus_\Delta J_\Delta^{ab}$$
where the sums are taken over all $n$-simplexes in $\Sigma$. The map in question is a $KK$-equivalence on each summand $J_\Delta$, equivariant for the action of the stabilizer group for $\Delta$ by 2.11. To obtain the inverse in $KK^G(I_n^{ab}/I_{n+1}^{ab}, I_n/I_{n+1})$ for the $KK^G$-element defined by the canonical map $I_n/I_{n+1} \to I_n^{ab}/I_{n+1}^{ab}$ we can take the direct sum (with the obvious $G$-action) of the Kasparov bimodules representing the inverse in $KK(J_\Delta^{ab}, J_\Delta)$ for the element defined by the quotient map $J_\Delta \to J_\Delta^{ab}$. $\square$

**Theorem 2.13** *Let $\Sigma$ be a finite-dimensional locally finite simplicial complex and $G$ a group acting simplicially on $\Sigma$. The canonical quotient map $\mathcal{C}_\Sigma \to \mathcal{C}_\Sigma^{ab}$ is a $KK^G$-equivalence.*

*Proof.* The skeleton filtration $(I_k)$ for $\mathcal{C}_\Sigma$ is finite. Since all quotient maps with respect to these ideals admit equivariant completely positive splittings by 2.10, we can use the exactness properties of $KK^G$ to deduce the assertion from repeated application of the five lemma (the spectral sequence associated with the filtration) from 2.12. $\square$

# 3    Noncommutative flag complexes.

In our applications to the Baum-Connes conjecture below we will be concerned with simplicial complexes that have a special property. We say, that a simplicial complex $\Sigma$ is full, if it is determined by its 1-simplexes through the following condition:

a finite subset $F \neq \emptyset$ of $V_\Sigma$ belongs to $\Sigma$ if $\{s,t\} \in \Sigma$ for all $\{s,t\} \subset F$.

Such complexes are sometimes called flag complexes in the literature. We will use both terminologies in the sequel. Note that the barycentric subdivision of any simplicial complex is full.

For a full locally finite complex we introduce a new noncommutative algebra $\mathcal{C}_\Sigma^{flag}$ as the universal $C^*$-algebra with generators $h_s, s \in V_\Sigma$ satisfying the relations $h_s \geq 0$ for all $s$, $h_s h_t = 0$ for $\{s,t\} \notin \Sigma$ and $\sum_{s \in V_\Sigma} h_s h_t = h_t, t \in V_\Sigma$. As before, the sum in the last condition is necessarily finite and this condition may be abbreviated formally to $\sum_{s \in V_\Sigma} h_s = 1$. The motivation for this definition is the fact that $C^*$-algebras associated in this way with flag complexes appear completely canonically in connection with the "finite" $K$-theory introduced in section 4.

Note that this $C^*$-algebra can also be defined for a simplicial complex that is not full. It will then coincide with the $C^*$-algebra $\mathcal{C}_{\Sigma^{flag}}^{flag}$ associated with the saturation $\Sigma^{flag}$ of $\Sigma$. (We define $\Sigma^{flag}$ to be the unique full simplicial complex with the same vertex set and the same edges as $\Sigma$. A subset $F$ of vertices belongs to $\Sigma^{flag}$ if $\{s,t\} \in \Sigma$ for all $s,t \in F$.)

For a flag complex, there are canonical surjective maps

$$\mathcal{C}_\Sigma^{flag} \longrightarrow \mathcal{C}_\Sigma \longrightarrow \mathcal{C}_\Sigma^{ab}$$

Unfortunately, without the skeleton filtration, the $C^*$-algebra $\mathcal{C}_\Sigma^{flag}$ does not, in general, retain the topological information of $\Sigma$. We are now going to determine the $K$-theory for a typical example of a noncommutative flag complex. This example is given by a complex that models a simplicial sphere.

Consider the universal $C^*$-algebra $C^*(x_0, x_1, \ldots, x_n)$ with self-adjoint generators $x_i$ satisfying the relation that $\sum_i x_i^2 = 1$. The abelianization of this $C^*$-algebra is obviously isomorphic to the algebra of continuous functions on the $n$-sphere $S^n$.

Consider the positive and negative parts $(x_i)_+$ and $(x_i)_-$ if $x_i$. We have $x_i^2 = (x_i)_+^2 + (x_i)_-^2$ and $(x_i)_+^2 (x_i)_-^2 = 0$. If we put $h_{i+} = (x_i)_+^2$ and $h_{i-} = (x_i)_-^2$ we get $2n+2$ new generators for $C^*(x_0, x_1, \ldots, x_n)$ that are positive and satisfy

$$\sum h_{i+} + \sum h_{i-} = 1, \qquad h_{i+} h_{i-} = 0 \quad \forall i$$

Thus $C^*(x_0, x_1, \ldots, x_n)$ is exactly isomorphic to the algebra $\mathcal{C}_\Sigma^{flag}$ for the flag complex $\Sigma_{S^n}$ with vertices $\{0^+, 1^+, \ldots, n^+, 0^-, \ldots, n^-\}$ and the condition that exactly the edges $\{i^+, i^-\}$ do not belong to $\Sigma_{S^n}$. This complex manifestly is a simplicial model for the $n$-sphere.

We will denote the $C^*$-algebra $\mathcal{C}^{flag}_{\Sigma_{S^n}}$ by $S^{nc}_n$.

For any vertex $s$ in a flag complex there is a natural evaluation map $\mathcal{C}^{flag}_\Sigma \to \mathbb{C}$ mapping the generator $h_s$ to 1 and all the other generators to 0. We will now consider various homomorphisms from $S^{nc}_n$ to $S^{nc}_n$ and from $S^{nc}_n$ to $M_2(S^{nc}_n)$.

Denote by $\kappa_i$ the homomorphism from $S^{nc}_n$ to $S^{nc}_n$ that maps the generator $h_{i+}$ to 1 and all the other generators to 0. This is simply the composition of the evaluation map for the vertex $i^+$ with the natural inclusion map $\mathbb{C} \to S^{nc}_n$.

We will further consider two homomorphisms $\alpha$ and $\beta$ from $S^{nc}_n$ to $M_2(S^{nc}_n)$ of the form

$$\alpha = \begin{pmatrix} \mathrm{id} & 0 \\ 0 & \alpha_1 \end{pmatrix} \qquad \beta = \begin{pmatrix} \beta_0 & 0 \\ 0 & \beta_1 \end{pmatrix}$$

where id is the identity homomorphism of $S^{nc}_n$ and $\alpha_1, \beta_0, \beta_1$ are the homomorphisms $S^{nc}_n \to S^{nc}_n$ determined by

$$\alpha_1(h_{i+}) = h_{i+} + h_{i-},\ \alpha_1(h_{i-}) = 0 \text{ for } i = 0, 1,$$

$$\alpha_1(h_s) = h_s \text{ for all other generators } h_s$$

$$\beta_0(h_{0+}) = h_{0+} + h_{0-},\ \beta_0(h_{0-}) = 0, \quad \beta_0(h_s)) = h_s \text{ for all other generators } h_s$$

$$\beta_1(h_{1+}) = h_{1+} + h_{1-},\ \beta_1(h_{1-}) = 0, \quad \beta_1(h_s)) = h_s \text{ for all other generators } h_s$$

**Lemma 3.1** *We have the following homotopies:*

*(a) $\alpha$ is homotopic to $\beta$.*

*(b) $\alpha_1$ is homotopic to $\kappa_0$.*

*(c) $\beta_0$ is homotopic to $\kappa_0$ and $\beta_1$ is homotopic to $\kappa_1$.*

*Proof.* (a) is proved by rotating $\beta(h_{0+})$ and $\beta(h_{0-})$ to $\alpha(h_{0+})$ and $\alpha(h_{0-})$, i.e. we use the homomorphisms $\varphi_t$ from $S^{nc}_n$ to $M_2(S^{nc}_n)$ mapping $h^+_0, h^-_0$ to $R_t \beta(h^+_0) R^*_t$, $R_t \beta(h^-_0) R^*_t$ and all the other generators $h_s$ to $\beta(h_s)$, where $R_t$ are rotation matrices, $t \in [0, \pi/2]$. We have $\varphi_0 = \beta$ and $\varphi_{\pi/2} = \alpha$.
(b) A homotopy between $\alpha_1$ and $\kappa_0$ is obtained by considering the homomorphisms mapping $h_{0+}$ to $t(h_{0+} + h_{0-}) + (1-t)1$, $h_{0-}$ to 0, $h_{1+}$ to $t(h_{1+} + h_{1-})$, $h_{1-}$ to 0 and all the other generators $h_s$ to $th_s$, $t \in [0, 1]$.
(c) A homotopy for $\beta_0$ is given by the homomorphisms mapping $h_{0+}$ to $t(h_{0+} + h_{0-}) + (1-t)1$, $h_{0-}$ to 0 and all the other generators $h_s$ to $th_s$. The homotopy for $\beta_1$ is then obtained by exchanging the role of 0 and 1. □

From this we can easily deduce the $K$-theory of $S^{nc}_n$.

**Proposition 3.2** *The evaluation map $S^{nc}_n \to \mathbb{C}$ at the vertex $1^+$ (and thus at any vertex) induces an isomorphism $K_*(S^{nc}_n) \longrightarrow K_*(\mathbb{C})$.*

*Proof.* We have to show that $K_*(\kappa_1)$ is an isomorphism. This follows from the fact that
$$K_*(\mathrm{id}) + K_*(\kappa_0) = K_*(\kappa_0) + K_*(\kappa_1)$$
which is a consequence of the preceding lemma. $\square$

In order to recover the topological information in $S_n^{nc}$ we have to use its natural filtration. Let $(I_k^{nc})$ be the natural skeleton filtration, i.e. $I_k^{nc}$ is the ideal in $S_n^{nc}$ generated by products containing at least $k+1$ different generators.

One can show that, for $n=1$, the natural maps $S_1^{nc}/I_2^{nc} \to \mathcal{C}(S^1)$ and $S_1^{nc}/I_3^{nc} \to \mathcal{C}(S^1)$ induce isomorphisms in $K$-theory. Moreover one has $K_0(I_2^{nc}) = K_0(I_3^{nc}) = \mathbb{Z}$ and $K_1(I_2^{nc}) = K_1(I_3^{nc}) = 0$.

For $n=2$, one can show that the natural map $S_2^{nc}/I_3^{nc} \to \mathcal{C}(S^2)$ induces an isomorphism in $K$-theory. However already here we have somewhat unexpectedly, that $K_0(I_3^{nc}/I_4^{nc}) \cong \mathbb{Z}^3$.

For higher $n$, the situation becomes much more complicated.

To end this section, we want to emphasize that the arguments used to determine the $K$-theory of $S_n^{nc}$ in this section are far from being equivariant for natural group actions on $S_n^{nc}$.

# 4 Applications to the Baum-Connes conjecture

Let $\Gamma$ be a discrete group. We denote by $\lambda_s$, $s \in \Gamma$ the unitary operators in $\mathcal{L}(\ell^2\Gamma)$ given by left translation $\lambda_s \xi_t = \xi_{st}$ on the canonical orthormal basis $\xi_s$ in $\ell^2\Gamma$. As usual, $C^*_{\mathrm{red}}\Gamma$ then denotes the $C^*$-algebra generated by the $\lambda_s$, $s \in \Gamma$. There is a natural coproduct $\delta_\Gamma : C^*_{\mathrm{red}}\Gamma \to C^*_{\mathrm{red}}\Gamma \otimes C^*_{\mathrm{red}}\Gamma$ given by $\delta_\Gamma(\lambda_s) = \lambda_s \otimes \lambda_s$ (all $C^*$-algebra tensor products in this paper will be minimal tensor products).

We are going to examine the $K$-theory of the reduced group $C^*$-algebra $C^*_{\mathrm{red}}\Gamma$ or, more generally, the $K$-theory of a reduced crossed product $A \rtimes_r \Gamma$ for any $C^*$-algebra $A$ carrying an action of $\Gamma$.

Now, such crossed products carry, as an additional structure, the dual action of $\hat{\Gamma}$, i.e. an action of the Hopf algebra $C^*_{\mathrm{red}}\Gamma$ (or a coaction of $\Gamma$).

An action of $\hat{\Gamma}$ on a $C^*$-algebra $B$ is, by definition, a nondegenerate (i.e. the hereditary subalgebra generated by the image is everything) homomorphism

$$\delta_B : B \to B \otimes C^*_{\mathrm{red}}\Gamma$$

satisfying the coassociativity condition $(\delta_B \otimes 1)\delta_B(x) = (1 \otimes \delta_\Gamma)\delta_B(x)$. Following [1] we say that $B$ is a $\hat{\Gamma}$-algebra if, in addition, $\delta_B$ is injective.

Any $\hat{\Gamma}$-algebra carries a natural $\Gamma$-grading, see [1]. A $C^*$-algebra $B$ is $\Gamma$-graded if for each $s \in \Gamma$, there is a closed linear subspace $B_s$ of $B$ such that

- $B = \overline{\bigoplus_{s \in \Gamma} B_s}$
- $B_s B_t \subset B_{st}$ for all $s, t \in \Gamma$

- $B_s^* = B_{s^{-1}}$

For a $\hat{\Gamma}$-algebra $B$ the projection onto the subspace $B_s$ is obtained by taking the coefficient at $s$ of $\delta_B(x)$. For $B = A \rtimes_r \Gamma$ we obviously can take $B_s = Au_s$, $u_s$ the unitary implementing $s \in \Gamma$.

Let $B$ be a $\hat{\Gamma}$-algebra. A $\hat{\Gamma}$-equivariant Hilbert-$B$-module is a Hilbert-$B$-module $E$ with a (coassociative nondegenerate) coaction $\delta_E : E \to E \otimes B$ satisfying

$$\delta_E(\xi)\delta_B(a) = \delta_E(\xi b) \qquad \delta_B((\xi \,|\, \eta)) = (\delta_E(\xi) \,|\, \delta(\eta))$$

for $\xi, \eta \in E$, $b \in B$ and the $B$-valued scalar product $(\,\cdot\,|\,\cdot\,)$, see [1], Definition 2.2 (this definition also makes sense if the coaction on $B$ is not injective). Using the coaction we can define subspaces $E_s$ of degree $s$ in $E$.

We now use the $\Gamma$-grading to define groups that approximate the $K$-theory of algebras with an action of $\hat{\Gamma}$. To define the finite and strongly finite $K$-theory of $B$ we consider triples $(p, \bar{p}, E)$ where $E$ is a $\hat{\Gamma}$-equivariant Hilbert $B$-module and $p, \bar{p}$ are projections in $\mathcal{L}(E)$ such that $p - \bar{p} \in \mathcal{K}(E)$. Here, as usual, $\mathcal{L}(E)$ and $\mathcal{K}(E)$ denote the algebras of bounded and compact operators on $E$, respectively. We denote by $\langle p, \bar{p}, E \rangle$ the equivalence class, for unitary $\hat{\Gamma}$-equivariant equivalence of such a triple. For any $\hat{\Gamma}$-equivariant Hilbert $B$-module $E$ we say that an element $x \in \mathcal{L}(E)$ has degree $s$ if $xE_t \subset E_{st}$ for all $t \in \Gamma$. This definition makes $\mathcal{K}(E)$ into a $\Gamma$-graded $C^*$-algebra. The algebra $\mathcal{L}(E)$ is not $\Gamma$-graded, since linear combinations of homogeneous elements with a fixed degree are not necessarily dense. However, for any $x \in \mathcal{L}(E)$, we can define its component $x_s$ of degree $s$ by $x_s \xi = P_{E_{st}} x \xi$ for $\xi \in E_t$ and $P_{E_{st}}$ the projection onto $E_{st}$. Let $x$ be an element and $A$ be a subalgebra of $\mathcal{L}(E)$. We define their $\Gamma$-support as

$$\operatorname{supp}_\Gamma x = \{s \in \Gamma \mid x_s \neq 0\}$$

$$\operatorname{supp}_\Gamma A = \{s \in \Gamma \mid A_x \neq 0\}$$

Let $F \subset \Gamma$ be a finite subset. Set

$$C_0^F(B) = \{\langle p, \bar{p}, E \rangle \,|\, \operatorname{supp}_\Gamma p,\ \operatorname{supp}_\Gamma \bar{p} \subset F\}$$

$$C_0^{fin}(B) = \{\langle p, \bar{p}, E \rangle \,|\, \operatorname{supp}_\Gamma p,\ \operatorname{supp}_\Gamma \bar{p} \text{ are finite}\}$$

$$C_0^{sF}(B) = \{\langle p, \bar{p}, E \rangle \,|\, \operatorname{supp}_\Gamma C^*(\{p_s | s \in \Gamma\}),\ \operatorname{supp}_\Gamma(C^*(\{\bar{p}_s | s \in \Gamma\})) \subset F\}$$

$$C_0^{sfin}(B) = \{\langle p, \bar{p}, E \rangle \,|\, C^*(\{p_s\}), C^*(\{\bar{p}_s\}) \text{ have finite } \Gamma - \text{support}\}$$

A homotopy between classes $\langle p_0, \bar{p}_0, E_0 \rangle$ and $\langle p_1, \bar{p}_1, E_1 \rangle$ in $C_0^F(B)$, $C_0^{\text{fin}}(B)$, $C_0^{sF}(B)$, $C_0^{\text{sfin}}(B)$, respectively, is as usually a class $\langle p, \bar{p}, E[0,1]\rangle$ in $C_0^F(B[0,1])$, $C_0^{\text{fin}}(B[0,1])$, $C_0^{sF}(B[0,1])$, $C_0^{\text{sfin}}(B[0,1])$, respectively, whose restrictions to 0 and $1 \in [0,1]$ are $\langle p_0, \bar{p}_0, E_0\rangle$ and $\langle p_1, \bar{p}_1, E_1\rangle$ (here, as usually, $B[0,1]$ denotes the algebra of $B$-valued continuous functions on the unit interval $[0,1]$). We denote by $[p_0, \bar{p}_0, E_0]$ the homotopy class of $\langle p_0, \bar{p}_0, E_0\rangle$.

We can now define for any $\hat{\Gamma}$-algebra $B$

**Definition 4.1**
$$K_0^F(B) = \{[p, \bar{p}, E] \mid \langle p, \bar{p}, E\rangle \in C_0^F(B)\}$$
$$K_0^{\text{fin}}(B) = \{[p, \bar{p}, E] \mid \langle p, \bar{p}, E\rangle \in C_0^{\text{fin}}(B)\}$$
$$K_0^{sF}(B) = \{[p, \bar{p}, E] \mid \langle p, \bar{p}, E\rangle \in C_0^{sF}(B)\}$$
$$K_0^{\text{sfin}}(B) = \{[p, \bar{p}, E] \mid \langle p, \bar{p}, E\rangle \in C_0^{\text{sfin}}(B)\}$$

With direct sum as addition $K_0^F(B), K_0^{\text{fin}}(B), K_0^{sF}(B)$ and $K_0^{\text{sfin}}(B)$ are abelian groups. Moreover, we obviously have

$$K_0^{\text{fin}}(B) = \varinjlim_F K_0^F(B)$$

$$K_0^{\text{sfin}}(B) = \varinjlim_F K_0^{sF}(B)$$

where the inductive limits are taken over all finite subsets $F$ of $\Gamma$.

There are natural maps $K_0^{\text{sfin}}B \to K_0^{\text{fin}}B \to K_0 B$. We may similarly define $K_1^{\text{sfin}}$ and $K_1^{\text{fin}}$ starting from triples $(p, T, E)$ where $E$ is a $\hat{\Gamma}$-equivariant Hilbert $C^*$-module, $p \in \mathcal{L}(E)$ is a projection with support condition using $F$ and $T \in \mathcal{L}(E)$ is an operator such that $T = T^*, T^2 = 1, T = T_1$ (i.e $\text{supp}_\Gamma T \subset \{1\}$) and $Tp - pT \in \mathcal{K}(E)$. Equivalently, we could define $K_1^{\text{sfin}}$ and $K_1^{\text{fin}}$ by taking suspensions.

Let now $E$ be a $\hat{\Gamma}$-equivariant Hilbert $B$-module and $p \in \mathcal{L}(E)$ a projection with finite $\Gamma$-support. The relations $p = p^*$ and $p^2 = p$, combined with the identity $p = \sum_{s \in \Gamma} p_s$ give, respectively:

$$(P_1) \quad p_s^* = p_{s^{-1}}$$
$$(P_2) \quad p_t = \sum_{s \in \Gamma} p_s p_{s^{-1}t}$$

Note that, in particular, $p_1^2 = \sum_{s \in \Gamma} p_s p_s^*$ ($1 \in \Gamma$ the neutral element). This shows that $p_1$ is a strictly positive element in the $C^*$-algebra generated by the $p_s, s \in \Gamma$.

**Definition 4.2** *Let $F$ be a finite subset in $\Gamma$. We denote by $P_F$ the universal $C^*$-algebra generated by elements $p_s, s \in F$, satisfying the relations $(P_1)$ and $(P_2)$ (where $p_s$ is understood to be zero for $s \notin F$). We denote by $P_F^s$ the universal $C^*$-algebra with generators $p_s, s \in F$ satisfying the relations $(P_1)$ and $(P_2)$ and the additional relation $p_{s_1} p_{s_2} \ldots p_{s_n} = 0$ whenever $s_1 s_2 \ldots s_n \notin F$.*

Note that for any projection $p$ as above, the relations $(P_1)$ and $(P_2)$ imply that $\text{supp}_\Gamma p = (\text{supp}_\Gamma p)^{-1}$ and $\text{supp}_\Gamma p \neq \emptyset \Rightarrow 1 \in \text{supp}_\Gamma p$. Therefore, in the definition of $P_F$ and $P_F^s$, we might just as well assume that $F = F^{-1}$ and $1 \in F$.

There is a homomorphism $\delta : P_F \to P_F \otimes C_{\text{red}}^* \Gamma$ defined by $\delta(p_s) = p_s \otimes \lambda_s$ ($\lambda_s \in C_{\text{red}}^* \Gamma$ the unitary defined by the left regular representation $\lambda$). Therefore, on $P_F$ we can define a natural $\Gamma$-grading (obtained by projecting $\delta(x)$ onto the coefficient of $\lambda_s$) such that $\deg p_s = s$. Also, $P_F$ contains a canonical projection, namely $p = \sum p_s$.

The same comments apply to $P_F^s$.

The algebra $P_F^s$ is a sum of finitely many closed subspaces $V_t$, $t \in F$, where $V_t$ is the closure of linear combinations of products $p_{s_1} p_{s_2} \ldots p_{s_k}$ with $s_1 s_2 \ldots s_k = t$. Each $V_t$ is mapped isometrically under the coaction $\delta : P_F^s \to P_F^s \otimes C_{red}^* \Gamma$ to $V_t \otimes \lambda_t$. Thus the coaction is injective and $P_F^s$ is a $\hat{\Gamma}$-algebra. This is not necessarily the case for $P_F$. We therefore define

**Definition 4.3** *We denote by $P_F^{\text{red}}$ the $C^*$-algebra isomorphic to the image of $P_F$ under the map $\delta : P_F \to P_F \otimes C_{\text{red}}^* \Gamma$, i.e. the $C^*$-algebra generated by the elements $p_s \otimes \lambda_s$, $s \in F$.*

$P_F^{\text{red}}$ is a $\hat{\Gamma}$-algebra. If $B$ is a $\hat{\Gamma}$-algebra, i.e. if the coaction $\delta : B \to B \otimes C_{\text{red}}^* \Gamma$ is injective, then any homomorphism $P_F \to B$ that respects the $\Gamma$-grading factors through $P_F^{\text{red}}$.

Let $E$ again be a $\hat{\Gamma}$-equivariant Hilbert module over the $\hat{\Gamma}$-algebra $B$. By the very definition of $P_F$ there is a bijection between:

- projections $e$ in $\mathcal{L}(E)$ with $\text{supp}_\Gamma e \subset F$
- homomorphisms $\varphi : P_F \to \mathcal{L}(E)$ respecting the $\Gamma$-grading

given by $\varphi \mapsto e_\varphi = \sum \varphi(p_s)$ and $e \mapsto \varphi_e$ where $\varphi_e$ is defined by $\varphi_e(p_s) = e_s$. If the coaction on $B$ (and thus on $E$) is injective, then these objects are moreover in bijection with

- homomorphisms $\varphi : P_F^{\text{red}} \to \mathcal{L}(E)$ respecting the $\Gamma$-grading

By definition of $P_F^s$ there is, similarly, a canonical bijection between:

- projections $e$ in $\mathcal{L}(E)$ with $\text{supp}_\Gamma(C^*(\{e_s | s \in \Gamma\})) \subset F$
- homomorphisms $\varphi : P_F^s \to \mathcal{L}(E)$ respecting the $\Gamma$-grading

Note also that a homomorphism from $P_F$ or $P_F^s$ into a $\hat{\Gamma}$-algebra $B$ is $\hat{\Gamma}$-equivariant if and only if it respects the $\hat{\Gamma}$-grading. We obtain, basically from the definition of the equivariant $KK$-groups, [1]

**Proposition 4.4** *Let $B$ be a $\hat{\Gamma}$-algebra and $F$ a finite subset of $\Gamma$. We have natural isomorphisms*

$$K_*^F(B) \cong KK_*^{\hat{\Gamma}}(P_F, B)$$

$$K_*^{fin}(B) \cong \varinjlim_F KK_*^{\hat{\Gamma}}(P_F, B)$$

$$K_*^{sF}(B) \cong KK_*^{\hat{\Gamma}}(P_F^s, B)$$

$$K_*^{sfin}(B) \cong \varinjlim_F KK_*^{\hat{\Gamma}}(P_F^s, B)$$

*The algebra $P_F$ can be replaced by $P_F^{\text{red}}$ in the first two isomorphisms.*

*Proof.* The proof follows simply from the fact that $KK_0^{\hat{\Gamma}}(B_1, B_2)$ can be defined as the set of homotopy classes of triples $(\varphi, \bar{\varphi}, E)$ where $E$ is a $\hat{\Gamma}$-equivariant Hilbert $B_2$-module and $\varphi, \bar{\varphi} : B_1 \to \mathcal{L}(E)$ are $\hat{\Gamma}$-equivariant homomorphisms, such that $\varphi(x) - \bar{\varphi}(x) \in \mathcal{K}(E)$ for all $x \in B_1$ (the operator $F$ appearing in the definition of $KK^{\hat{\Gamma}}$ in [1] can here be assumed to be $\hat{\Gamma}$-invariant). In the odd case, $KK_1^{\hat{\Gamma}}(B_1, B_2)$ can similarly be defined from triples $(\varphi, T, E)$ where $E, \varphi$ are as above and $T \in \mathcal{L}(E)$ satisfies $T = T^*, T^2 = 1, \text{supp}_\Gamma T \subset \{1\}$ and $T\varphi(x) - \varphi(x)T \in \mathcal{K}(E)$ for all $x \in B_1$. These descriptions of $KK_*^\Gamma$ show that the claimed isomorphisms are simply a translation of the definitions of $K_*^F, K_*^{fin}, K_*^{sF}$ and $K_*^{sfin}$. $\square$

We are now in a position to apply the Baaj-Skandalis duality. Recall that, for any $\hat{\Gamma}$-algebra $A$, one can define a reduced crossed product $A \rtimes_r \hat{\Gamma}$. The precise definition is irrelevant for our purposes. We only need the following two basic results, [1],[12]:

$$KK_*^{\hat{\Gamma}}(B_1, B_2) \cong KK_*^\Gamma(B_1 \rtimes_r \hat{\Gamma}, B_2 \rtimes_r \hat{\Gamma})$$

$$A \rtimes_r \Gamma \rtimes_r \hat{\Gamma} \cong A \otimes \mathcal{K}(\ell^2 \Gamma)$$

which hold for all $\hat{\Gamma}$-algebras $B_1$, $B_2$ and for each $\Gamma$-algebra $A$. We obtain for $B = A \rtimes_r \Gamma$

**Proposition 4.5** *Let $A$ be a $\Gamma$-algebra. Then*

$$K_*^{fin}(A \rtimes_r \Gamma) \cong \varinjlim_F KK_*^\Gamma(P_F^{\text{red}} \rtimes_r \hat{\Gamma}, A)$$

$$K_*^{sfin}(A \rtimes_r \Gamma) \cong \varinjlim_F KK_*^\Gamma(P_F^s \rtimes_r \hat{\Gamma}, A)$$

*Proof.* Combining the Baaj-Skandalis results with 4.5 we obtain isomorphisms

$$K_*^F(A\rtimes_r\Gamma) \cong KK_*^{\hat{\Gamma}}(P_F^{\text{red}}, A\rtimes_r\Gamma) \cong KK_*^\Gamma(P_F^{\text{red}}\rtimes_r\hat{\Gamma}, A\rtimes_r\Gamma\rtimes_r\hat{\Gamma}) \cong KK_*^\Gamma(P_F^{\text{red}}\rtimes_r\hat{\Gamma}, A)$$

and analogous isomorphisms for $K_*^{sF}$. $\square$

We are now going to determine $P_F^{\text{red}} \rtimes_r \hat{\Gamma}$ and $P_F^s \rtimes_r \hat{\Gamma}$ (up to Morita equivalence). For this we define

- $\mathcal{E}_F$ is the universal $C^*$-algebra with generators $h_s$, $s \in \Gamma$ satisfying the relations $h_s \geq 0, h_s h_t = 0$ if $s^{-1}t \notin F$ and $\sum_s h_s h_t = h_t$

- $\mathcal{E}_F^s$ is the quotient of $\mathcal{E}_F$ by the additional relation

$$h_s h_{a_1} h_{a_2} \ldots h_{a_n} h_t = 0$$

  for all $a_1, \ldots, a_n \in \Gamma$, if $s^{-1}t \notin F$

These algebras are exactly noncommutative simplicial complex algebras, associated with a natural full simplicial complex, in the sense of sections 2 and 3.

**Proposition 4.6** *Let $\Sigma_F$ be the full locally finite simplicial complex with vertex set $\Gamma$ and given by*

$$\Sigma_F = \{\{s_0, \ldots, s_n\} \subset \Gamma \mid s_i^{-1}s_j \in F, \ \forall i,j = 0, \ldots, n\}$$

*Then $\mathcal{E}_F^s \cong \mathcal{C}_{\Sigma_F}$ and $\mathcal{E}_F \cong \mathcal{C}_{\Sigma_F}^{flag}$.*

*Proof.* Obvious from the definitions. □

Let $\mathcal{E}_F^{ab}$ denote the abelianization of $\mathcal{E}_F$, i.e. the quotient of $\mathcal{E}_F$ by the additional relation $h_s h_t = h_t h_s, s, t \in \Gamma$. We have natural maps

$$\mathcal{E}_F \longrightarrow \mathcal{E}_F^s \longrightarrow \mathcal{E}_F^{ab}$$

We say that two $C^*$-algebras $A$ and $B$ with an action of $\hat{\Gamma}$ are $\hat{\Gamma}$-equivariantly Morita equivalent, if there is a $C^*$-algebra $E$ carrying an action of $\hat{\Gamma}$ with a direct sum decomposition into closed $\hat{\Gamma}$-invariant subspaces $A, X, Y$ and $B$ that multiply like in a $2 \times 2$-matrix decomposition

$$\begin{pmatrix} A & X \\ Y & B \end{pmatrix}$$

i.e. $A$ and $B$ are subalgebras of $E$ such that $AB = BA = 0$, $AX = X$, $BY = Y$, $XY = A$, $YX = B$.
Any $\hat{\Gamma}$-equivariant Morita equivalence between two algebras $A$ and $B$ with an action of $\hat{\Gamma}$ induces an invertible element in $KK^{\hat{\Gamma}}(A, B)$, [1] 5.9.

**Theorem 4.7** *Let $F$ be a finite subset of $\Gamma$ (with $F = F^{-1}$ and $1 \in F$). The $C^*$-algebra $P_F$ is $\hat{\Gamma}$-equivariantly Morita equivalent to the full crossed product $\mathcal{E}_F \rtimes \Gamma$ while $P_F^{\text{red}}$ is isomorphic to the reduced crossed product $\mathcal{E}_F \rtimes_r \Gamma$. The $C^*$-algebra $P_F^s$ is $\hat{\Gamma}$-equivariantly Morita equivalent to $\mathcal{E}_F^s \rtimes \Gamma \cong \mathcal{E}_F^s \rtimes_r \Gamma$.*

*Proof.* There is a natural $\Gamma$-grading respecting homomorphism $\alpha$ from $P_F$ onto a $\hat{\Gamma}$-invariant $C^*$-subalgebra of $\mathcal{E}_F \rtimes \Gamma$ mapping $p_s \in P_F$ to $h_1^{\frac{1}{2}} u_s h_1^{\frac{1}{2}}$ ($u_s$ the unitary implementing the automorphism of $\mathcal{E}_F$ given by $s$). The identity

$$h_1 h_{s_1} h_{s_2} \ldots h_{s_n} u_t h_1 = \alpha(p_1^{\frac{1}{2}} p_{s_1} p_{s_1^{-1} s_2} \ldots p_{s_{n-1}^{-1} s_n} p_{s_n^{-1} t} p_1^{\frac{1}{2}})$$

shows that $\alpha(P_F) = \overline{h_1 (\mathcal{E} \rtimes \Gamma) h_1}$. It follows that $\alpha(P_F)$ is a full hereditary subalgebra of $\mathcal{E}_F \rtimes \Gamma$.
To prove the first assertion, it will therefore suffice to show that $\alpha$ is injective (in the definition of Morita equivalence we can then take $X = \overline{h_1(\mathcal{E}_F \rtimes \Gamma)}$ and $Y = \overline{(\mathcal{E}_F \rtimes \Gamma) h_1}$).
Let $P_F$ be faithfully represented as a subalgebra of $\mathcal{L}(H)$ and let $\varrho_s, \lambda_s, s \in \Gamma$ denote the unitary operators in $\mathcal{L}(l^2 \Gamma)$ defined by the right and left regular representation of Γ, respectively.

Consider the operator

$$\hat{p} = \sum_{s \in F} p_s \otimes \varrho_s$$

in $\mathcal{L}(H \otimes l^2\Gamma)$. The relations satisfied by the $p_s$ immediately imply that $\hat{p}$ is a projection (up to replacing $\rho_s$ by $\lambda_s$, we have $\hat{p} = \delta_\Gamma(p)$).

Moreover, $\hat{p}$ is invariant under the action of $\Gamma$ by conjugation by $1 \otimes \lambda_t$, $t \in \Gamma$.

Let $e_{st}$ for $s,t \in \Gamma$ denote the matrix units in $\mathcal{K}(l^2\Gamma)$ (mapping the $t$-th element in the canonical orthonormal basis for $l^2\Gamma$ to the $s$-th element) and put

$$f_s = \hat{p}(1 \otimes e_{ss})\hat{p}$$

Then $f_s f_t = \hat{p}(1 \otimes e_{ss})\hat{p}(1 \otimes e_{tt})\hat{p} = \hat{p}(p_{s^{-1}t} \otimes e_{st})\hat{p}$.

Thus $f_s f_t = f_t f_s = 0$ for $s^{-1}t \notin F$ and moreover $\sum_{s \in \Gamma} f_s f_t = f_t$ (since $\sum_s (1 \otimes e_{ss})\hat{p}(1 \otimes e_{tt}) = \hat{p}(1 \otimes e_{tt})$).

Also, obviously $(1 \otimes \lambda_s) f_t (1 \otimes \lambda_s^{-1}) = f_{st}$.

Mapping the generators $h_s$ of $\mathcal{E}_F$ to $f_s$, we therefore obtain a $\Gamma$-covariant representation of $\mathcal{E}_F$ in $H \otimes l^2\Gamma$, whence a representation of $\mathcal{E}_F \rtimes \Gamma$.

In order to show that the homomorphism $\alpha$, defined at the beginning of the proof, is injective, it is enough to show that the homomorphism $\beta : P_F \to \mathcal{L}(H \otimes \ell^2\Gamma)$ mapping $p_s$ to $f_1^{\frac{1}{2}}(1 \otimes \lambda_s) f_1^{\frac{1}{2}}$ is injective (note that $\beta$ factors through $\alpha$).

This however follows from the identity

$$(1 \otimes e_{11})(f_1 \beta(p_{s_1} p_{s_2} \ldots p_{s_n}) f_1)(1 \otimes e_{11}) = (p_1 p_{s_1} p_{s_2} \ldots p_{s_n} p_1) \otimes e_{11}$$

for all $s_1, s_2, \ldots, s_n \in F$, which implies more generally that

$$(1 \otimes e_{11}) f_1 \, \beta(x) \, f_1 \, (1 \otimes e_{11}) = (p_1 x p_1) \otimes e_{11} \quad \text{for all } x \in P_F$$

and the fact that the positive map $P_F \ni x \mapsto p_1 x p_1 \in P_F$ is injective ($p_1$ is a strictly positive element in $P_F$).

To treat the case of $P_F^{\text{red}}$ one simply has to replace $p_s$ by $p_s \otimes \lambda_s$ everywhere.

The argument for $P_F^s$ is also a verbatim repetition of the preceding discussion. One checks that, in that case

$$f_s f_{a_1} \ldots f_{a_n} f_t = \hat{p} \left( p_{s^{-1}a_1} p_{a_1^{-1}a_2} \ldots p_{a_n^{-1}t} \otimes e_{st} \right) \hat{p}$$

so that the $f_s$ satisfy the relations of $\mathcal{E}_F^s$. □

**Corollary 4.8** *Let $F$ be a finite subset of $\Gamma$ (with $F = F^{-1}$ and $1 \in F$). The $C^*$-algebra $\mathcal{E}_F$ is $\Gamma$-equivariantly Morita equivalent to $P_F^{\text{red}} \rtimes_r \hat{\Gamma}$. The $C^*$-algebra $\mathcal{E}_F^s$ is $\Gamma$-equivariantly Morita equivalent to $P_F^s \rtimes_r \hat{\Gamma}$. We have*

$$K_*^{fin}(A \rtimes_r \Gamma) \cong \varinjlim_F KK_*^\Gamma(\mathcal{E}_F, A)$$

$$K_*^{sfin}(A \rtimes_r \Gamma) \cong \varinjlim_F KK_*^\Gamma(\mathcal{E}_F^s, A)$$

*Proof.* This follows from 4.7 by Baaj-Skandalis duality. □

Recall now the definition of the Baum-Connes map for a $C^*$-algebra $A$ with an action of $\Gamma$. The topological $K$-theory of $\Gamma$ with coefficients in $A$ is defined as $KK_*^\Gamma(\underline{E}_\Gamma, A)$ where $\underline{E}_\Gamma$ is the universal $\Gamma$-space with a proper $\Gamma$-action. A natural realization of $\underline{E}_\Gamma$ is

$$\underline{E}_\Gamma = \{\,\mu \,|\, \mu \text{ is a unit measure on } \Gamma \text{ with finite support}\,\}$$

The technical definition of $KK_*^\Gamma(\underline{E}_\Gamma, A)$ now is

$$KK_*^\Gamma(\underline{E}_\Gamma, A) = \varinjlim_F KK_*^\Gamma(\mathcal{C}_0(\underline{E}_F), A)$$

where

$$\underline{E}_F = \{\,\mu \,|\, \mu \text{ is a unit measure on } \Gamma \text{ with } \operatorname{supp}\mu \subset s_o F \text{ for some } s_0 \in \Gamma\,\}$$

and where the inductive limit is taken over all finite subsets $F \subset \Gamma$, see e.g.[6]. Note now that $\underline{E}_F$ is exactly the spectrum of $\mathcal{E}_F^{ab}$.

**Proposition 4.9** *The $C^*$-algebra $\mathcal{E}_F^{ab}$ is isomorphic to $\mathcal{C}_0(\underline{E}_F)$. We have*

$$KK_*^\Gamma(\underline{E}_\Gamma, A) \cong \varinjlim_F KK_*^\Gamma(\mathcal{E}_F^{ab}, A)$$

*The Baum-Connes map $\mu : KK_*^\Gamma(\underline{E}_\Gamma, A) \to K_*(A \rtimes_r \Gamma)$ is exactly the composition of the maps*

$$KK_*^\Gamma(\underline{E}_\Gamma, A) = \varinjlim_F KK_*^\Gamma(\mathcal{E}_F^{ab}, A) \to \varinjlim_F KK_*^\Gamma(\mathcal{E}_F, A) \cong K_*^{fin}(A \rtimes_r \Gamma) \to K_*(A \rtimes_r \Gamma)$$

*where the first arrow is induced by the natural map $\mathcal{E}_F \to \mathcal{E}_F^{ab}$.*

*Proof.* The Baum-Connes map is constructed by choosing a so called cut-off function $f$ on $\underline{E}_F$. A cut-off function is a positive function with compact support such that the sum over all translates $s(f)$, $s \in \Gamma$ is the constant function 1 and such that $fs(f) \neq 0$ only for finitely many $s$. One then considers the projection $e = \sum_{s \in \Gamma} f^{\frac{1}{2}} u_s f^{\frac{1}{2}} \in \mathcal{C}_0(\underline{E}_F) \rtimes_r \Gamma$. The map $\mu$ then maps an element $\alpha \in KK_*^\Gamma(\underline{E}_F, A)$ to the Kasparov product $[e] \cdot [\alpha \rtimes_r \Gamma] \in K_*(A \rtimes_r \Gamma)$ where $[\alpha \rtimes_r \Gamma]$ is the element of $KK_*(\mathcal{C}_0(\underline{E}_F) \rtimes_r \Gamma, A \rtimes_r \Gamma)$ obtained from $\alpha$ by "descent".

Now, the generator $h_1 \in \mathcal{E}_F^{ab} \cong \mathcal{C}_0(\underline{E}_F)$ is precisely a possible cut-off function. Under the natural map $P_F^{\text{red}} \subset \mathcal{E}_F \rtimes_r \Gamma \to \mathcal{C}_0(\underline{E}_F) \rtimes_r \Gamma$ the universal projection $p$ is mapped to $e$. If $u$ and $v$ denote the natural maps $K_*(P_F^{\text{red}}) \to K_*(\mathcal{E}_F^{ab} \rtimes_r \Gamma)$ and $KK_*(\mathcal{E}_F^{ab} \rtimes_r \Gamma, A \rtimes_r \Gamma) \to KK_*(P_F^{\text{red}}, A \rtimes_r \Gamma)$, then $u([p]) \cdot [\alpha \rtimes_r \Gamma] = [p] \cdot v([\alpha \rtimes_r \Gamma])$. However, $[p] \cdot v([\alpha \rtimes_r \Gamma])$ is exactly the image of $v([\alpha \rtimes_r \Gamma])$ under the natural map $KK_*(P_F^{\text{red}}, A \rtimes_r \Gamma) \to K_*(A \rtimes_r \Gamma)$. □

The commutative diagram

$$\begin{array}{ccc} KK^\Gamma_*(\mathcal{E}^s_F, A) & \longrightarrow & K^{sfin}_*(A \rtimes_r \Gamma) \\ \downarrow & & \downarrow \\ KK^\Gamma_*(\mathcal{E}_F, A) & \longrightarrow & K^{fin}_*(A \rtimes_r \Gamma) \end{array}$$

shows that the Baum-Connes map is also equal to the composition of the maps

$$KK^\Gamma_*(\underline{E}_\Gamma, A) \xrightarrow{\varphi} \varinjlim_F KK^\Gamma_*(\mathcal{E}^s_F, A) \cong K^{sfin}_*(A \rtimes_r \Gamma) \xrightarrow{\psi} K_*(A \rtimes_r \Gamma).$$

We will now show that $\varphi$ is always an isomorphism. Therefore the Baum-Connes map $\mu$ is an isomorphism if and only if the natural map

$$K^{sfin}_*(A \rtimes_r \Gamma) \xrightarrow{\psi} K_*(A \rtimes_r \Gamma)$$

is an isomorphism.

**Proposition 4.10** *Let $F$ be a finite subset of $\Gamma$. The canonical map $\mathcal{E}^s_F \to \mathcal{E}^{ab}_F$ is a $KK^\Gamma$-equivalence.*

*Proof.* This follows from 2.13, since the dimension of $\Sigma_F$ is bounded by $|F|-1$. □

We summarize the preceding discussion in the following theorem, which is our main result .

**Theorem 4.11** *Let $A$ be a $C^*$-algebra on which $\Gamma$ acts by automorphisms. The left hand side $KK^\Gamma_*(\underline{E}_\Gamma, A)$ of the Baum-Connes conjecture is canonically isomorphic to $K^{sfin}_*(A \rtimes_r \Gamma)$. The Baum-Connes map $\mu : KK^\Gamma_*(\underline{E}_\Gamma, A) \to K_*(A \rtimes_r \Gamma)$ corresponds, under this isomorphism, to the natural map $K^{sfin}_*(A \rtimes_r \Gamma) \to K_*(A \rtimes_r \Gamma)$.*

We thus have the following factorizations for the Baum-Connes map:

$$KK^\Gamma_*(\underline{E}_\Gamma, A) \cong K^{sfin}_*(A \rtimes_r \Gamma) \longrightarrow K_*(A \rtimes_r \Gamma).$$

and

$$KK^\Gamma_*(\underline{E}_\Gamma, A) \longrightarrow K^{sfin}_*(A \rtimes_r \Gamma) \longrightarrow K^{fin}_*(A \rtimes_r \Gamma) \longrightarrow K_*(A \rtimes_r \Gamma)$$

In view of the discussion in section 3, it seems quite uncertain if the map

$$K^{sfin}_*(A \rtimes_r \Gamma) \longrightarrow K^{fin}_*(A \rtimes_r \Gamma)$$

could be an isomorphism in general. We mention however that, if $A$ is proper in the sense of [6], then it is nearly trivial that the map $K^{fin}_*(A \rtimes_r \Gamma) \to K_*(A \rtimes_r \Gamma)$ is an isomorphism. Thus, by the result in [6] combined with Theorem 4.11, the map $K^{sfin}_*(A \rtimes_r \Gamma) \to K^{fin}_*(A \rtimes_r \Gamma)$ has to be an isomorphism.
The connection between $K^{fin}_*(A \rtimes_r \Gamma)$ and $K^{sfin}_*(A \rtimes_r \Gamma)$ becomes however much closer if we take the skeleton filtration into account. Given a finite subset $F$ of $\Gamma$, let $(I_k)$ and $(I^s_k)$ be the skeleton filtrations in $\mathcal{E}_F$ and $\mathcal{E}^s_F$, respectively. We denote by $J_F$ the kernel of the natural map $\mathcal{E}_F \to \mathcal{E}^s_F$.

**Lemma 4.12** *Let $F$ be a finite subset of $\Gamma$ with $F = F^{-1}$ and $1 \in F$, and $F^n$ the set of all products of $n$ elements in $F$. Let moreover denote $h_s$, $s \in \Gamma$, the generators of $\mathcal{E}_F$ and $\mathcal{E}_F^s$.*

(a) *If $s^{-1}t \notin F^n$ for $s, t$ in $\Gamma$, then for all $t_1, \ldots, t_{n-1}$ in $\Gamma$, the product*

$$h_s h_{t_1} \ldots h_{t_{n-1}} h_t$$

*is zero in $\mathcal{E}_F$.*

(b) *The natural map $\mathcal{E}_{F^n} \longrightarrow \mathcal{E}_F/I_{n+1}$ maps $J_{F^n}$ to zero.*

(c) *The natural map $\mathcal{E}_{F^n}/I_{n+1} \longrightarrow \mathcal{E}_F/I_{n+1}$ factors as*

$$\mathcal{E}_{F^n}/I_{n+1} \longrightarrow \mathcal{E}_{F^n}^s/I_{n+1}^s \longrightarrow \mathcal{E}_F/I_{n+1}$$

*Proof.* Obvious. □

**Proposition 4.13** *Let $A$ be a $C^*$-algebra with an action of $\Gamma$. For each fixed $n$, we have an isomorphism*

$$\varinjlim_F KK_*^\Gamma(\mathcal{E}_F/I_n, A) \cong \varinjlim_F KK_*^\Gamma(\mathcal{E}_F^s/I_n^s, A) \tag{1}$$

*Proof.* This follows from 4.12 (c). □

If we take the inductive limit over $n$ in the right hand side of equation 1 we obtain $K_*^{sfin}(A \rtimes_r \Gamma)$.